\documentclass[10pt,a4paper]{article}
\usepackage{amssymb,latexsym,amsmath,amsfonts,amsthm}
\makeatletter
\newtheorem{thm}{Theorem}[section]
\newtheorem{lem}[thm]{Lemma}
\newtheorem{prop}[thm]{Proposition}

\theoremstyle{remark}

\newtheorem{ex}{\textsc{Example}}[section]
\newcommand{\RR}{\mathbb R}

\newcommand{\Hspace}{\mathcal H}

\def\dint{-\kern-1.1em\int}

\catcode`@=11
\newbox\tr@tto
\setbox\tr@tto=\hbox{{\count0=0\dimen0=-,9pt\dimen1=1,1pt\loop\ifnum
    \count0<11 \advance \count0 by1 \vrule width.51pt height\dimen1
    depth\dimen0\kern-0.17pt\advance\dimen0 by-0.05pt\advance\dimen1
    by0.1pt\repeat \loop\ifnum\count0<21\advance \count0 by1 \vrule
    width.6pt height\dimen1 depth\dimen0\kern-0.2pt \advance\dimen0
    by-0.1pt\advance\dimen1 by 0.05pt\repeat}}
\def\medint{\displaystyle\copy\tr@tto\kern-10.4pt\int}
\catcode`@=12

\makeatletter \protected@xdef\@thefnmark{\relax}
\def\AMSclass#1{\@footnotetext{\hspace{-4mm}--- \it{Mathematics Subject Classification} \rm{(2000)}:
 #1.}}
\def\keywords#1{\@footnotetext{\hspace{-4mm}--- \it{Keywords:}
 #1.}}
\begin{document}
\title{A sharp weighted Wirtinger inequality and some related functional spaces}
\author{Raffaella Giova\\
\small{Dipartimento di Statistica e  Matematica per la Ricerca Economica}\\
\small{Universit\`a di Napoli ``Parthenope",}
\small{Via Medina, 40 - 80133 Napoli, Italy}\\
\small{Fax: +39-081\,5474904; E-mail: raffaella.giova@uniparthenope.it}
\and Tonia Ricciardi\thanks{
Corresponding author. Partially supported by Regione Campania L.R.~5/2006, by
GNAMPA-INdAM,
and by a PRIN-COFIN project
of the Italian Ministry of University and Research.}\\
\small{Dipartimento di Matematica e Applicazioni ``R. Caccioppoli"}\\
\small{Universit\`a di Napoli Federico II,}
\small{Via Cintia - 80126 Napoli, Italy}\\
\small{Fax: +39-081\,675665; E-mail: tonia.ricciardi@unina.it}\\
}
\date{March 10, 2008}
\maketitle
\begin{abstract}
We consider the generalized
Wirtinger inequality
\[
\left( \int_{0}^{T} a |u|^q \right)^{1/q} \le C \biggm(\int_{0}^{T}
a^{1-p} |u'|^{p}\biggm)^{1/p},
\]
with $p,q>1$, $T>0$, $a\in L^1[0,T]$, $a\ge0$, $a\not\equiv0$ and
where $u$ is a $T$-periodic function satisfying the constraint
\[
\int_{0}^{T} a |u|^{q-2}u =0.
\]
We provide the best constant $C>0$ as well as all extremals.
Furthermore, we characterize the natural functional space where
the inequality is defined.
\end{abstract}
\begin{description}
\item {\textsc{Key Words:}} Weighted Wirtinger inequality, best constant,
weighted Sobolev space, generalized trigonometric functions.
\item {\textsc{MSC 2000 Subject Classification:}} 26D15.
\end{description}
%\vskip.6truecm
%\subjclass{26D15}
%%%%%%%%%%%%%%%%%%%%%%%%%%%%%%%%%%%%%%%%%%%%%%%%%%%%%%%%%%%%%%%%%%%%%%%%%%%%%%%%%%%%%%%%%%%%%%%
%%%%%%%%%%%%%%%%%%%%%%%%%%%%%%%%%%%%%%%%%%%%%%%%%%%%%%%%%%%%%%%%%%%%%%%%%%%%%%%%%%%%%%%%%%%%%%%
%%%%%%%%%%%%%%%%%%%%%%%%%%%%%%%%%%%%%%%%%%%%%%%%%%%%%%%%%%%%%%%%%%%%%%%%%%%%%%%%%%%%%%%%%%%%%%%
%%%%%%%%%%%%%%%%%%%%%%%%%%%%%%%%%%%%%%%%%%%%%%%%%%%%%%%%%%%%%%%%%%%%%%%%%%%%%%%%%%%%%%%%%%%%%%%
\section{Introduction and main results}
Wirtinger type inequalities are of interest in various areas of
analysis and mathematical physics, including the Wulff theorem
\cite{DP}, quasiconformal mapping theory \cite{IS}, $p$-Laplacian
systems \cite{MM}. In view of their applications, they received a
considerable attention in recent years. See, e.g.,
\cite{CD,DGS,Gi1,Li,Ri1} and the references therein.
\par
In reference~\cite{Ri1}, the following weighted Wirtinger inequality is considered:
\begin{equation}
\label{wirt22} \int_0^{2\pi}au^2\le C\int_0^{2\pi}a^{-1}u'^2,
\end{equation}
where $a\ge0$ and $u$ is a
$2\pi$-periodic function satisfying $\int_0^{2\pi}au=0$. By a
technique introduced in \cite{PS}, the sharp value of the constant
$C>0$ is a key ingredient in \cite{Ri1}, which is used in order to obtain a sharp H\"older estimate
for two-dimensional, divergence-form, elliptic equations with unit
determinant coefficient matrix. Such equations are closely related
to quasiconformal mappings and Beltrami equations, see
\cite{IS,Ri2}. In this context, the coefficient $a$, which is
related to the coefficient matrix of the elliptic equation,
naturally satisfies the assumption $a,a^{-1}\in L^\infty[0,2\pi]$. By the
diffeomorphism $y=2\pi(\int_0^{2\pi}a)^{-1}\int_0^xa(t)\,dt$,
inequality~\eqref{wirt22} reduces to the case $a\equiv1$, which is
the standard Wirtinger inequality. Thus, we obtain
$C=\left((2\pi)^{-1}\int_0^{2\pi}a\right)^2$. Moreover, the
extremals are given by $u(x)=\alpha\sin\left(\widetilde
a^{-1}\int_0^xa\,dt+\delta\right)$ for some $\alpha\neq0$ and
$\delta\in\RR$, where $\widetilde a=(2\pi)^{-1}\int_0^{2\pi}a\,dt$.
Such a value of $C$ suggests that inequality~\eqref{wirt22} should
hold true in the more general case $a\in L^1(0,2\pi)$. This fact was
established, among other results, in \cite{Gi1}. Furthermore, in \cite{Gi2} the best
constant $C>0$ in the more general inequality
\begin{equation}
\label{wirtpp} 
\int_0^{2\pi}a|u|^p\le C\int_0^{2\pi}a^{1-p}|u'|^p
\end{equation}
is given.
Here, $a\in L^1(0,2\pi)$, $a\ge0$, $a\not\equiv0$ and $u$ is a
$2\pi$-periodic function satisfying $\int_0^{2\pi}a|u|^{p-2}u=0$.
We note that the proofs in \cite{Gi1,Gi2} are based on an approximation argument
involving truncations, which
does not allow to characterize the extremals.
\par
In this note we consider the following generalized Wirtinger inequality
\begin{equation}
\label{wirtpq} 
\left(\int_0^{T}a|u|^q\right)^{1/q}\le
C\left(\int_0^{T}a^{1-p}|u'|^p\right)^{1/p}
\end{equation}
where $p,q>1$, $T>0$, $a\ge0$, $a\not\equiv0$, $a\in L^1[0,T]$ and where
$u$ is a $T$-periodic function satisfying $\int_0^Ta|u|^{q-2}u=0$.
When $p=q$, \eqref{wirtpq} reduces to \eqref{wirtpp}.
We note that, although the underlying reason for which inequality~\eqref{wirtpq}
holds is a rescaling argument, the rigorous analysis is not obvious
under our general assumptions on $a$. Indeed, $a$ may vanish on a
set of positive measure. In fact, one of our problems is to identify a
natural space $\mathcal X$ where inequality \eqref{wirtpq} is
defined. In the space $\mathcal X$ we can characterize all extremals
in terms of generalized trigonometric functions \cite{Hu,Li,LP,MM}.
Then, we show that $\mathcal X$ includes the natural weighted
Sobolev spaces defined in the usual way, by approximation by smooth
functions or by distributional derivatives. We show that for some
particular choices of $a$, such weighted Sobolev spaces may be
strictly included in $\mathcal X$.
\par
More precisely, we define:
\[
\mathcal X=\left\{u:[0,T]\to\RR:\ \exists U\in W_{\mathrm{per}}^{1,p}(\RR)\
\mbox{such that}\ u(x)=U\left(\frac{T}{\int_0^Ta}\int_0^xa\right)\right\},
\]
where $W_{\mathrm{per}}^{1,p}(\RR)=\{u\in W^{1,p}(\RR):\ u\mathrm{\ is\ }T\mathrm{-periodic}\}$. 
With this notation, we have:
\begin{thm}
\label{thm:main} Let $p,q>1$ and $T>0$. Let $a \in L^1[0,T]$,
$a\ge0$, $a\not\equiv0$. Then, the following Wirtinger inequality
holds:
\begin{equation}
\label{mainineq} \biggm(\int_{0}^{T} a |u|^q \biggm)^{1/q} \le
\widetilde a^{1/p^*+1/q}\,C(p,q)\biggm(\int_{0}^{T} a^{1-p}
|u'|^{p}\biggm)^{1/p}
\end{equation}
for every $u \in\mathcal X$ such that
\[
\int_{0}^{T} a |u|^{q-2}u=0,
\]
where $\widetilde a=T^{-1}\int_{0}^{T}a$ and
\begin{equation}
\label{Capq} C(p,q)=\left[2\left(\frac{1}{p^*}\right)^{1/q}
\left(\frac{1} {q} \right)^{1/{p^*}}\left( \frac 2
{p^*+q}\right)^{1/p - 1/q}B \left( \frac {1}{p^*} , \frac {1} {q}
\right)\right]^{-1}.
\end{equation}
\end{thm}
Throughout this note, for every $p\ge1$ we denote by $p^*=p/(p-1)$ 
the conjugate exponent of $p$.
Moreover, $B(\alpha,\beta)$ denotes the Beta function
\[
B(\alpha,\beta)=\int_0^1 t^{\alpha-1}(1-t)^{\beta-1}\,dt=B(\beta,\alpha)
\]
for all $\alpha,\beta>1$.
\par
In Section~\ref{sec:spaces} we show that $\mathcal X$ is a Banach space
and we compare it to the weighted Sobolev spaces, as defined in the
usual ways.
More precisely, setting
\[
\mathcal A=\left\{u\in C^1(\RR):\ u\ \mathrm{is\ }T\mathrm{-periodic\ and\ }
\int_{0}^{T}a^{1-p}|u'|^p<+\infty\right\},
\]
we show that $\|\cdot\|$ defined by
\[
\|u\|=\left(\int_0^Ta|u|^q\right)^{1/q}+\left(\int_0^Ta^{1-p}|u'|^p\right)^{1/p}
\]
is a norm on $\mathcal A$ (recall that we allow $a$ to vanish on a set of positive measure).
We define $\Hspace$ as the closure of $\mathcal A$ with respect to $\|\cdot\|$.
We show that $\Hspace\subset\mathcal X$.
Finally, we compare $\mathcal X$ with the space
\[
\mathcal W=\left\{u\in L^1[0,T]:\ \int_0^Ta^{1-p}|u'|^p<+\infty
\ \mathrm{and\ }u(0)=u(T)\right\},
\]
which was considered in \cite{Gi1,Gi2}.
We show that $\mathcal W=\mathcal X$, so that the assumption
$u\in L^1[0,T]$, which is needed in order to define the distributional derivative of $u$,
but does not seem natural in \eqref{wirtpq}, is actually not restrictive.
We also show that for particular choices of $a$, we may have
$\Hspace\neq\mathcal W$, unlike what happens in the usual Sobolev spaces,
see \cite{MS}.
It follows in particular that Theorem~\ref{thm:main} holds for all functions
belonging to the traditional weighted Sobolev spaces.
\par
In order to characterize the extremals for
\eqref{mainineq}, we use generalized trigonometric functions, which
we now briefly define. See \cite{Hu,Li,LP,MM} for more details. Let
$p,q>1$. The function ${\arcsin}_{pq}:[0,1]\to\RR$ is defined by
\[
{\arcsin}_{pq}(\sigma)= \int_0^{\sigma}\frac{ds}{(1-s^p)^{1/q^*}}.
\]
Then, we have
\[
{\arcsin}_{pq}(1)= \frac 1 p B\left(\frac 1 p,\frac 1 q
\right)=:\frac {\pi_{pq}} 2.
\]
The function ${\arcsin}_{pq}:[0,1] \rightarrow [0,\frac{\pi_{pq}} 2]$
is strictly increasing and its inverse function is
denoted by $\sin_{pq}$.
The function $\sin_{pq}$ is extended as an odd function
to the interval $[-\pi_{pq},\pi_{pq}]$ by setting
${\sin}_{pq}(t)={\sin}_{pq}(\pi_{pq}-t)$ in $[\pi_{pq}/2,\pi_{pq}]$,
$\sin_{pq}(t)=-\sin_{pq}(-t)$ in $[-\pi_{pq},0]$,
and to the whole real axis as
a $2\pi_{pq}-$periodic function.
The function $\xi(t)=\sin_{qp^*}(\pi_{qp^*}t)$ is the unique solution of the
initial value problem:
\[
(|u'|^{p-2}u')'+\frac q{p^*}(|u|^{q-2}u)=0, \qquad u(0)=0,\quad
u'(0)=1
\]
and it satisfies:
\begin{align}
\label{pqbasic} \frac{\|\xi'\|_p}{\|\xi\|_q}
=&\inf\Big\{\frac{\|u'\|_p}{\|u\|_q}:\
u\in W^{1,p}(\RR)\setminus\{0\},\ u\mathrm{\ is\ 2-periodic}\\ 
\nonumber
&\qquad\qquad 
\qquad\qquad\mathrm{and\ }\int_{-1}^{1}|u|^{q-2}u=0\Big\}\\
\nonumber
=&C^{-1}(p,q),
\end{align}
where $C(p,q)$ is the constant defined in \eqref{Capq}. Moreover,
any minimizer for \eqref{pqbasic} is of the form
$u(t)=\alpha\,\xi(t+\delta)$ for some $\alpha\neq0$ and
$\delta\in\RR$. It should be mentioned that the existence of the
minimum in \eqref{pqbasic} and the explicit value of $C(p,q)$ were
obtained in \cite{DGS} in a more general setting, and further
generalized in \cite{CD}.
\par
At this point, we can state the sharpness of Theorem~\ref{thm:main}.
\begin{thm}
\label{thm:sharp}
Let $u\in\mathcal X$. Then $u$ satisfies inequality~\eqref{mainineq}
with the equal sign if and only if $u$ is of the form
\begin{equation}
\label{extremal}
u(x)=\alpha\,\sin_{qp^*}\left(\pi_{qp^*}
\,\frac{T}{\int_{0}^{T}a}\int_0^x a\,dt+\delta\right),
\end{equation}
for some $\alpha \in \RR\setminus \{0 \}$ and for some $\delta\in\RR$.
\end{thm}
The remaining part of this note is devoted to the proofs
of Theorem~\ref{thm:main} and Theorem~\ref{thm:sharp}.
%%%%%%%%%%%%%%%%%%%%%%%%%%%%%%%%%%%%%%%%%%%%%%%%%%%%%%%%%%%%%%%%%%%%%%%%%%%%%%%%%%%%%%%%%
%%%%%%%%%%%%%%%%%%%%%%%%%%%%%%%%%%%%%%%%%%%%%%%%%%%%%%%%%%%%%%%%%%%%%%%%%%%%%%%%%%%%%%%%%
%
%%%%%%%%%%%%%%%%%%%%%%%%%%%%%%%%%%%%%%%%%%%%%%%%%%%%%%%%%%%%%%%%%%%%%%%%%%%%%%%%%%%%%%%%%
%%%%%%%%%%%%%%%%%%%%%%%%%%%%%%%%%%%%%%%%%%%%%%%%%%%%%%%%%%%%%%%%%%%%%%%%%%%%%%%%%%%%%%%%%
\section{Some weighted Sobolev spaces}
\label{sec:spaces}
\subsection{The space $\mathcal W$}
The following space was considered in \cite{Gi1,Gi2}:
\[
\mathcal W=\left\{u\in L^1[0,T]:\ \int_0^Ta^{1-p}|u'|^p<+\infty
\ \mathrm{and\ }u(0)=u(T)\right\},
\]
where $u'$ denotes the distributional derivative of $u$.
We first check that $\mathcal W$ is well-defined.
\begin{lem}
\label{lem:Hoelder}
Let $a\in L^1[0,T]$ and let $f:[0,T]\to\RR$ be a measurable function.
Then,
\begin{equation*}
\label{Hoelder}
\int_{0}^{T} |f|\le \biggm(\int_{0}^{T}|f|^p
a^{1-p}\biggm)^{1/p}\biggm(\int_{0}^{T} a\biggm)^{(p-1)/p}.
\end{equation*}
\end{lem}
\begin{proof}
The proof follows by the H\"older inequality.
\end{proof}
In view of Lemma~\ref{lem:Hoelder}, we have that $u'\in L^1[0,T]$
for all $u\in\mathcal W$. It follows that $u$ is absolutely
continuous. In particular, the periodicity condition $u(0)=u(T)$ is
well-defined. When $p=q$, it was shown in \cite{Gi1,Gi2} that
Theorem~\ref{thm:main} holds for all $u\in\mathcal W$ by an
approximation argument. A natural question is whether one can define
a suitable functional space for inequality \eqref{mainineq} which
does not require $u\in L^1[0,T]$. Indeed, this condition is only
used to define the distributional derivative of $u$, and does not
seem natural for \eqref{mainineq}.
\par
In what follows, we consider some other natural functional spaces.
%%%%%%%%%%%%%%%%%%%%%%%%%%%%%%%%%%%%%%%%%%%%%%%%%%%%%%%%%%%%%%%%%%%%%%%%%%%%%%%%%%%%%%%%%
%%%%%%%%%%%%%%%%%%%%%%%%%%%%%%%%%%%%%%%%%%%%%%%%%%%%%%%%%%%%%%%%%%%%%%%%%%%%%%%%%%%%%%%%%
\subsection{The space $\Hspace$}
For $k\ge0$, let
$C_{\mathrm{per}}^k(\mathbb R)$ denote the space of $C^k$ functions
defined on $\RR$ which are $T$-periodic. Let $s>1$, $b\in L^1[0,T]$,
$b \ge 0$ and $b\not\equiv0$. Let $|\cdot|_{b,s}$ be the seminorm
defined on $C_{\mathrm{per}}(\RR)$ by
\[
|u|_{b,s}=\biggm(\int_{0}^{T}  b|u|^s \biggm)^{1/s}.
\]
Let
\[
\mathcal A=\left\{u\in C_{\mathrm{per}}^1(\RR):\
\int_{0}^{T}a^{1-p}|u'|^p<+\infty\right\}.
\]
For each $u \in\mathcal A$ we define
\[
\|u \|_{b,s}=|u|_{b,s}+\biggm(\int_{0}^{T} a^{1-p}
|u'|^{p}\biggm)^{1/p}.
\]
\begin{lem}
\label{lem:norm} 
Let $a,b\in L^1[0,T]$, $a,b\ge0$, $a,b\not\equiv0$.
Then, $\|\cdot \|_{b,s}$ is a norm on $\mathcal A$.
\begin{proof}
We need only prove that $\|u\|_{b,s} =0$ implies $u=0$ for every
$u\in \mathcal A$. To this end, we observe that $\| u\|_{b,s}=0$ if and only if
$\int_{0}^{T}  b|u|^s=0$ and $\int_{0}^{T} a^{1-p} |u'|^{p}=0$.
Then, since $a^{1-p}>0$ a.e.\ in $[0,T]$, it results $u'=0$ a.e.\
and therefore, since $u'$ is continuous, $u' \equiv 0$. It follows
that $u\equiv c\equiv\mathrm{const}$. From $\int_{0}^{T} b>0$, we
conclude that $\int_{0}^{T} b|u|^s=c^s \int_{0}^{T} b=0$ and
therefore $c=0$.
\end{proof}
\end{lem}
We denote by $\Hspace_{b,s}$ the closure of
$\mathcal A$ with respect to the norm $\|\cdot\|_{b,s}$.
\begin{prop}
\label{prop:spaces} Let $s,{\widetilde{s}}>1$, $b,{\widetilde{b}}\in
L^1(0,T)$, $b,{\widetilde{b}} \ge0$ and $b,
{\widetilde{b}}\not\equiv0$. Then
\[
\Hspace_{b,s}=\Hspace_{{\widetilde{b}},{\widetilde{s}}}
\subset L^\infty[0,T]
\]
as sets of functions.
\end{prop}
\begin{proof}
Let $u\in\Hspace_{b,s}$.
We note that
$\int_{0}^{T}|u'|<\infty$. Indeed,
by Lemma~\ref{lem:Hoelder}, we have:
\begin{align}
\label{wprimeL1}
\int_{0}^{T} |u'|\le&\biggm(\int_{0}^{T}|u'|^p
a^{1-p}\biggm)^{1/p}\biggm(\int_{0}^{T} a\biggm)^{(p-1)/p}\\
\nonumber
 \le&
\|u\|_{b,s}\biggm(\int_{0}^{T} a\biggm)^{(p-1)/p}<+\infty.
\end{align}
Let $ u_n \in C^1_{\mathrm{per}}(\RR)$ be a Cauchy sequence in
$\Hspace_{b,s}$. Then, $\int_{0}^{T}  b|u_n-u_m|^s \rightarrow
0$ and $\int_{0}^{T} a^{1-p} |u'_n-u'_m|^{p} \rightarrow 0$, as $m,n
\rightarrow \infty$. Furthermore, by \eqref{wprimeL1} with $u=u_n-u_m$,
we have
\[
\int_{0}^{T}  |u'_n-u'_m|\le \biggm(\int_{0}^{T} a^{1-p}
|u'_n-u'_m|^{p}\biggm)^{1/p}\biggm(\int_{0}^{T} a\biggm)^{(p-1)/p}
\rightarrow 0.
\]
Let
\[
E=\{x\in [0,T]:\ b(x)>0\}.
\]
We observe that the seminorm $|\cdot|_{b,s}$ restricted  to the
functions defined on $E$ defines an ordinary weighted Lebesgue space $L_b^s(E)$.
Then ${u_n}|_E$ converges in $L_b^s(E)$ and there exists
a subsequence $u_{n_k}$ which converges a.e.\ in $E$. Let $x_0\in E$ be such that
$u_{n_k}(x_0)$ converges. By the fundamental theorem of calculus,
for all $y \in [0,T]$ we have
$u_{n_k}(y)=u_{n_k}(x_0)+\int_{x_0}^{y}u'_{n_k}$,
$u_{n_h}(y)=u_{n_h}(x_0)+\int_{x_0}^{y}u'_{n_h}$ and consequently
\begin{equation*}
u_{n_k}(y)-u_{n_h}(y)=u_{n_k}(x_0)-u_{n_h}(x_0)+\int_{x_0}^{y}(u'_{n_k}-u'_{n_h}).
\end{equation*}
It follows that
\[
\sup_{y\in [0,T]}|u_{n_k}(y)-u_{n_h}(y)| \le
|u_{n_k}(x_0)-u_{n_h}(x_0)|+\int_{0}^{T}|u'_{n_k}-u'_{n_h}|\rightarrow
0
\]
as $h,k\to\infty$.
Therefore $u_{n_k}$ is a Cauchy sequence in $L^{\infty}[0,T]$ and
there exists $v \in C_{\mathrm{per}}(\RR)$  such that
$u_{n_k}\rightarrow v$ in $L^{\infty}[0,T]$. We now prove that the
whole sequence $u_n$ converges to $v$ in $L^{\infty}[0,T]$. To this
end, let $u_{n_{h,1}}$ and $u_{n_{h,2}}$ be subsequences of $u_n$
such that
\[
u_{n_{h,1}} \rightarrow v_1 \qquad \quad u_{n_{h,2}} \rightarrow
v_2,
\]
with $v_1,v_2\in C[0,T]$. We have:
\begin{align*}
&\left(\int_{0}^{T}b|u_{n_{h,1}}-u_{n_{h,2}}|^s\right)^{1/s}
\ge\left(\int_{0}^{T}b|v_1-v_2|^s\right)^{1/s}-\\
&\qquad\qquad\left(\int_{0}^{T}b|u_{n_{h,1}}-v_1|^s\right)^{1/s}-
\left(\int_{0}^{T}  b|u_{n_{h,2}}-v_2|^s\right)^{1/s}.
\end{align*}
Since $\int_{0}^{T}b|u_{n_{h,1}}-u_{n_{h,2}}|^s \rightarrow 0$,
$\int_{0}^{T}b|u_{n_{h,1}}-v_1|^s \rightarrow 0$,
$\int_{0}^{T}b|u_{n_{h,2}}-v_2|^s \rightarrow 0$ we obtain
$\int_{0}^{T} b|v_1-v_2|^s\to0$.
That is, $\displaystyle{\int_{0}^{T} } b|v_1-v_2|^s=0$ and
consequently $v_1=v_2$ in $E$. We have proved that any convergent
subsequence of $u_n$ converges to $v$ in $L^\infty(E)$. We conclude
that $u_n\to v$ in $L^\infty(E)$. We fix $x_0\in E$. By the
fundamental theorem of calculus, for all $y\in[0,T]$ we have
$u_n(y)=u_n(x_0)+\int_{x_0}^{y}u'_n$,
$u_m(y)=u_m(x_0)+\int_{x_0}^{y}u'_m$ and therefore
\begin{equation*}
u_n(y)-u_m(y)=u_n(x_0)-u_m(x_0)+\int_{x_0}^{y}(u'_n-u'_m).
\end{equation*}
It follows that, for all $y\in[0,T]$:
\[
|u_n(y)-u_m(y)|\le
|u_n(x_0)-u_m(x_0)|+\int_{0}^{T}|u'_n-u'_m|\rightarrow 0.
\]
Equivalently,
\[
\|u_n-u_m\|_{\infty}=\sup_{y \in [0,T]}|u_n(y)-u_m(y)|\rightarrow0
\]
as $m,n\to\infty$.
We conclude that $u_n$ is a Cauchy sequence in $L^{\infty}[0,T]$ and
$u_n\rightarrow v$ in $L^{\infty}[0,T]$ for some continuous function
$v$. 
At this point
it is readily seen that $u_n$ is a Cauchy sequence in 
$\Hspace_{\widetilde b,\widetilde s}$. Indeed, we have:
\[
\displaystyle{\int_{0}^{T} } \widetilde{b}|u_n-u_m|^{\widetilde{s}}
\le \|u_n-u_m\|_{\infty}^{\widetilde{s}} \displaystyle{\int_{0}^{T}
}\widetilde{b} \rightarrow 0.
\]
We conclude that a Cauchy sequence in $\Hspace_{b,s}$
is also a Cauchy sequence in $\Hspace_{\widetilde b,\widetilde s}$. Therefore,
as sets, $\Hspace_{b,s}=\Hspace_{\widetilde b,\widetilde s}$.
\end{proof}
In view of Proposition~\ref{prop:spaces}, we set
\[
\Hspace=\Hspace_{b,s}
\]
for any $b\in L^1[0,T]$, $b\ge0$, $b\not\equiv0$ and for any $s>1$.
%%%%%%%%%%%%%%%%%%%%%%%%%%%%%%%%%%%%%%%%%%%%%%%%%%%%%%%%%%%%%%%%%%%%%%%%%%%%%%%%%%%%
%%%%%%%%%%%%%%%%%%%%%%%%%%%%%%%%%%%%%%%%%%%%%%%%%%%%%%%%%%%%%%%%%%%%%%%%%%%%%%%%%%%%
% Proofs
%%%%%%%%%%%%%%%%%%%%%%%%%%%%%%%%%%%%%%%%%%%%%%%%%%%%%%%%%%%%%%%%%%%%%%%%%%%%%%%%%%%%
%%%%%%%%%%%%%%%%%%%%%%%%%%%%%%%%%%%%%%%%%%%%%%%%%%%%%%%%%%%%%%%%%%%%%%%%%%%%%%%%%%%%
\section{Proof of Theorem~\ref{thm:main} and of Theorem~\ref{thm:sharp}}
\label{sec:proofs}
Let $a\in L^1[0,T]$, $a\ge0$, $a\not\equiv0$ and let $p>1$.
We consider $y:[0,T]\to[0,T]$ defined by
\[y(x)=\widetilde a^{-1}\int_{0}^xa(t)\,dt,
\]
where $\widetilde a=T^{-1}\int_0^Ta$.
The function $y$ is well-defined,
nondecreasing, absolutely continuous and differentiable a.e.
We denote by $W_{\mathrm{per}}^{1,p}$ the set of functions in $W^{1,p}(\RR)$
which are $T$-periodic.
For every $U\in W_{\mathrm{per}}^{1,p}$ we define $u(x)=(\Psi U)(x)=U(y(x))$. By
taking difference quotients, we see that $u$ is differentiable a.e.,
and
\[
u'(x)=U'(y(x))\frac{a(x)}{\widetilde a}
\]
for almost every $x\in[0,T]$.
We recall the following general version of the change of variables
formula, see, e.g., \cite{HS}, Theorem~9.7.5, p.245.
\begin{lem}[\cite{HS}]
\label{lem:HS}
Let $g$ be a nondecreasing absolutely continuous function on $[\alpha,\beta]$
and let $f$ be integrable on $[g(\alpha),g(\beta)]$.
Then, $(f\circ g)g'\in L^1[\alpha,\beta]$ and
\[
\int_\alpha^\beta(f\circ g)g'=\int_{g(\alpha)}^{g(\beta)}f.
\]
\end{lem}
We set
\[
\mathcal X=\Psi(W_{\mathrm{per}}^{1,p}).
\]
\begin{lem}
\label{lem:X}
The mapping $\Psi:W_{\mathrm{per}}^{1,p}\to\mathcal X$ is
an isomorphism of Banach spaces.
\end{lem}
\begin{proof}
In view of Lemma~\ref{lem:HS}, we have:
\begin{align}
\label{changeu}
\int_0^T a|u|^q\,dx=&\widetilde a\int_0^T|U|^q\,dy\\
\label{changeu'}
\int_0^T a^{1-p}|u'|^p\,dx=&\widetilde a^{1-p}\int_0^T|U'|^p\,dy.
\end{align}
We note that $\mathcal X$ is complete with respect to the norm
$\|\cdot\|_{a,q}$ defined in Section~\ref{sec:spaces}. Indeed, 
in view of \eqref{changeu}--\eqref{changeu'},
if $u_n=\Psi(U_n)\in\mathcal X$ is
a Cauchy sequence, then $U_n$ is a Cauchy sequence in $W_{\mathrm{per}}^{1,p}$. Then
$U_n\to U$ in $W_{\mathrm{per}}^{1,p}$ and $u_n\to u=\Psi(U)$. Now the
claim follows by the Open Mapping Theorem.
\end{proof}
In the next lemma we clarify the relations between $\Hspace,\mathcal W$
and $\mathcal X$.
\begin{lem}
\label{lem:subset}
There holds:
\[
\Hspace\subset\mathcal W=\mathcal X.
\]
\end{lem}
\begin{proof}
We first show that $\Hspace\subset\mathcal X$. 
Let $u\in\Hspace$. We may assume that $u\in C^1$.
In the degenerate
case where $\inf a=0$, the function $y(x)$ may have some ``flat
regions". Namely, there may be $x',x''\in I$, $x'<x''$ such that
$y(x')=y(x'')=y(x)$ for all $x\in[x',x'']$. In this case, $a=0$
a.e.\ in $[x',x'']$. Since $\int_0^Ta^{1-p}u'<+\infty$, 
we derive that $u=$const in $[x',x'']$ and in particular
$u(x')=u(x'')$. It follows that $U(y)=u(x(y))$ is well-defined and
continuous and moreover $u(x)=U(y(x))$. 
Furthermore, since $x(y)$ is monotone, it is differentiable
a.e. It follows that $U$ is differentiable a.e. By change of
variables, as in Lemma~\ref{lem:HS},
$U\in W_{\mathrm{per}}^{1,p}$ and $u=\Psi U$.
Hence, $u\in\mathcal X$.
\par
Now we show that $\mathcal W\subset\mathcal X$. Let $u\in\mathcal
W$. Let $x_1\in[0,T]$ be a jump discontinuity point for $x(y)$.
Then, there exists $x_2>x_1$ such that $y(x)=y(x_1)=y(x_2)$ for all
$x\in[x_1,x_2]$. In particular, $a=0$ for almost every
$x\in[x_1,x_2]$ and consequently $u'=0$ for almost every
$x\in[x_1,x_2]$. It follows that $u(x)=u(x_1)=u(x_2)$ for every
$x\in[x_1,x_2]$. Therefore, the function $U(y)=u(x(y))$ is
well-defined and continuous. Moreover, consideration of difference
quotients yields $U'(y)=u'(x(y))x'(y)$ for almost every $y\in[0,T]$.
We claim that the almost everywhere derivative $U'$ is the
distributional derivative of $U$. To this end, let $\varphi\in
C_{\mathrm{per}}^1(\RR)$. Since $u'$ is the distributional
derivative of $u$, we have
\begin{align*}
\int_0^TU'(y)\varphi(y)\,dy
=&\int_0^Tu'(x)\varphi(y(x))\,dx
=-\int_0^Tu(x)\varphi'(y(x))y'(x)\,dx\\
=&-\int_0^TU(y)\varphi'(y)\,dy.
\end{align*}
Hence, $u=\Psi U$ with $U\in W_{\mathrm{per}}^{1,p}$ and therefore 
$u\in\mathcal X$.
\par
Finally, we show that $\mathcal X\subset\mathcal W$. To this end,
let $U\in W_{\mathrm{per}}^{1,p}(\RR)$ and let $u(x)=U(y(x))$. Then,
$u$ is continuous and $T$-periodic. Moreover, by taking difference
quotients, we have that $u$ is differentiably a.e.\ in $[0,T]$ and
the pointwise derivative is given by
\begin{equation*}
\label{uprime}
u'(x)=U'(y(x))y'(x)
\end{equation*}
for a.e.\ $x\in[0,T]$.
We have to show that $u'$ is the distributional derivative
of $u$. Since $U'$ is the distributional derivative of $U$,
for any $y_1\in[0,T]$ we have $U(y_1)=U(0)+\int_0^{y_1}U'(y)\,dy$.
Let $x_1=\inf\{x\in[0,T]:\ y(x)=y_1\}$. By the change
of variables formula, Lemma~\ref{lem:HS}, we obtain
$u(x_1)=u(0)+\int_0^{x_1}u'(x)\,dx$.
Let $x_2=\sup\{x\in[0,T]:\ y(x)=y_1\}$. Then,
$y'(x)=0$ for all $x\in(x_1,x_2)$ and in view of \eqref{uprime},
we have $u(x)=u(0)+\int_0^{x}u'(t)\,dt$ for all $x\in[x_1,x_2]$.
Since $x(y)$ only admits jump discontinuities, we conclude that
$u(x)=u(0)+\int_0^{x}u'(t)\,dt$ for all $x\in[0,T]$.
It follows that $u'$ is indeed the distributional derivative of $u$.
In view of \eqref{changeu'}, we conclude that $u\in\mathcal W$.
\end{proof}
By the following example we see that for some particular choices of
$a$, 
the space $\Hspace$ may degenerate to the space of constant functions,
and in particular $\Hspace\neq\mathcal W$,
unlike
what happens in the usual Sobolev spaces, see \cite{MS}.
\begin{ex}
There exists $a\in L^1[0,T]$, $a\ge0$, $a\not\equiv0$
such that 
\[
\Hspace=\{c\}_{c\in\RR}\neq\mathcal W.
\]
\end{ex}
Indeed, let $\mathcal C\subset[0,T]$ be a Cantor
set such that $|\mathcal C|=T/2=|[0,T]\setminus\mathcal C|$. Let
$a=\chi_\mathcal C$, the characteristic function of $\mathcal C$. We
claim  that $\mathcal X=\{c\}_{c\in\RR}$. To see this, recall that
$\mathcal C$ and $[0,T]\setminus\mathcal C$ are dense in $[0,T]$. Let
$u\in\mathcal A$, where $\mathcal A$ is defined in
Section~\ref{sec:spaces}. Since $a^{1-p}=+\infty$ on
$[0,T]\setminus\mathcal C$, we have $u'=0$ on $[0,T]\setminus\mathcal C$. By
continuity, $u'=0$ on $[0,T]$. It follows that $u$ is constant. Now let
$u\in\mathcal X$. Then, there exists $u_n=c_n\in\mathcal A$ such
that $c_n\to u$ with the respect to the norm $\|\cdot\|_{a,p}$. It
follows that $c_n$ is bounded, $c_n\to c\in\RR$ and $u=c$. We
conclude that for this choice of $a$, $\Hspace$ is the space of
constant functions.
\par
Finally, we provide the proofs of our main results.
\begin{proof}[Proof of Theorem~\ref{thm:main}]
Let $u\in\mathcal X$. Then $u(x)=U(y(x))$
for some $U\in W_{\mathrm{per}}^{1,p}(\RR)$.
The functions $u,U$ satisfy the identities \eqref{changeu}--\eqref{changeu'}
and moreover:
\begin{align*}
\int_{0}^{T} a |u|^{q-2}u=\widetilde{a}\int_{0}^{T}|U|^{q-2}U=0.
\end{align*}
Therefore, using \eqref{pqbasic} we conclude the
proof.
\end{proof}
\begin{proof}[Proof of Theorem~\ref{thm:sharp}]
Uniqueness of the extremals in $W_{\mathrm{per}}^{1,p}(\RR)$ implies uniqueness
of the extremals of the form~\eqref{extremal} in $\mathcal X$.
\end{proof}
%%%%%%%%%%%%%%%%%%%%%%%%%%%%%%%%%%%%%%%%%%%%%%%%%%%%%%%%%%%%%%%%%%%%%%%%%%%%%%%%%%%%
%%%%%%%%%%%%%%%%%%%%%%%%%%%%%%%%%%%%%%%%%%%%%%%%%%%%%%%%%%%%%%%%%%%%%%%%%%%%%%%%%%%%
%%%%%%%%%%%%%%%%%%%%%%%%%%%%%%%%%%%%%%%%%%%%%%%%%%%%%%%%%%%%%%%%%%%%%%%%%%%%%%%%%%%%
%%%%%%%%%%%%%%%%%%%%%%%%%%%%%%%%%%%%%%%%%%%%%%%%%%%%%%%%%%%%%%%%%%%%%%%%%%%%%%%%%%%%
\section*{Acknowledgement}
T.R.\ thanks Professor Jean Mawhin for providing her with
references \cite{Hu} and \cite{MM}.
%%%%%%%%%%%%%%%%%%%%%%%%%%%%%%%%%%%%%%%%%%%%%%%%%%%%%%%%%%%%%%%%%%%%%%%%%%%%%%%%%%%%
%%%%%%%%%%%%%%%%%%%%%%%%%%%%%%%%%%%%%%%%%%%%%%%%%%%%%%%%%%%%%%%%%%%%%%%%%%%%%%%%%%%%
%%%%%%%%%%%%%%%%%%%%%%%%%%%%%%%%%%%%%%%%%%%%%%%%%%%%%%%%%%%%%%%%%%%%%%%%%%%%%%%%%%%%
%%%%%%%%%%%%%%%%%%%%%%%%%%%%%%%%%%%%%%%%%%%%%%%%%%%%%%%%%%%%%%%%%%%%%%%%%%%%%%%%%%%%

\end{document}